\documentclass{svmult}

\usepackage{amsmath}
\usepackage{amssymb}
\usepackage[english]{babel}
\usepackage{tikz}
\usetikzlibrary{intersections, shapes}
\usetikzlibrary{arrows}
\usepackage{subcaption}
\usepackage{booktabs,multirow}
\usepackage{bigdelim}
\usepackage[hidelinks]{hyperref}
\usepackage{array}
\usepackage{pgfplots}
\usepackage[leftcaption]{sidecap}
\sidecaptionvpos{figure}{t}
\usepackage{graphicx,bm,xcolor}

\theoremstyle{definition} %

\begin{document}
\title*{
Towards a IETI-DP solver on non-matching multi-patch domains}
\titlerunning{}
\author{Rainer Schneckenleitner \inst{*}
	Stefan Takacs 
}
\institute{
	Rainer Schneckenleitner \at Institute of Computational Mathematics, Johannes Kepler University Linz, Altenberger Straße 69, 4040 Linz, Austria
	\email{schneckenleitner@numa.uni-linz.ac.at}, \\ 
	Corresponding author
	\and Stefan Takacs \at RICAM, Austrian Academy of Sciences, Altenberger Straße 69, 4040 Linz, Austria
	\email{stefan.takacs@ricam.oeaw.ac.at}
}
\maketitle
%
%
\abstract{
Recently, the authors have proposed and analyzed isogeometric tearing and interconnecting (IETI-DP) solvers for multi-patch discretizations in
Isogeometric Analysis. Conforming and discontinuous Galerkin
settings have been considered. In both cases, we have assumed that the
interfaces between the patches consist of whole edges.
In this paper, we present a generalization that allows us to drop this
requirement. This means that the patches can meet in T-junctions,
which increases the flexibility of the geometric model significantly.
We use vertex-based primal degrees of freedom.
For the T-junctions, we propose to follow the idea of ``fat vertices''.
}
%
%
\section{Introduction}
Isogeometric Analysis (IgA), see~\cite{HughesCottrellBazilevs:2005},
is a method for discretizing partial differential equations (PDEs).
The goal of its development has been to enhance the interface between computer-aided design (CAD) and simulation. Current state-of-the-art CAD tools use B-splines and NURBS for the representation of the
computational domain. 
In IgA, the same kind of bases 
is also utilized to discretize the PDEs.
Complex domains for real-world applications are usually the union
of many patches, parametrized with individual geometry functions (multi-patch IgA).
We focus on non-overlapping patches.

If the grids are not conforming and/or the interfaces between the patches do not consist of whole edges then discontinuous Galerkin (dG) methods are the discretization techniques of choice. A well studied representative is the symmetric interior discontinuous Galerkin (SIPG) method, cf.~\cite{Arnold:1982}. It has already been adapted and analyzed in IgA, cf.~\cite{LangerMantzaflaris:2015,LangerToulopoulos:2015,SchneckenleitnerTakacs:2020} and others. 
An obvious choice to solve discretized PDEs on domains with many non-overlapping patches are tearing and interconnecting methods. The variant we are interested in is the dual-primal approach,
see~\cite{FarhatLesoinneLeTallecPiersonRixen:2001a} for FETI-DP
and \cite{KleissPechsteinJuttlerTomar:2012,Hofer:2016a,HoferLanger:2016a}
for its extension to IgA, which
is called accordingly dual-primal isogeometric tearing and interconnecting method (IETI-DP). In~\cite{SchneckenleitnerTakacs:2019,SchneckenleitnerTakacs:2020}, the authors have presented a $p$- and $h$-robust convergence analysis. The authors have assumed that the interfaces consist of whole edges. If the vertices are chosen as primal degrees of freedom, it was shown that
the condition number of the preconditioned Schur complement system
is, under proper assumptions, bounded by
\begin{equation}\label{estimate}
	C\, p \; 
	\left(1+\log p+\max_{k=1,\ldots,K} \log\frac{H_{k}}{h_{k}}\right)^2,
\end{equation}
where $p$ is the spline degree, $h_k$ is the grid size on patch $\Omega^{(k)}$
and $H_k$ is the diameter of $\Omega^{(k)}$ and $C>0$ is a constant independent
of these quantities.
In this paper, we construct a new IETI-DP method that can deal
with interfaces that do not consist of whole edges.
This means that the patches can meet in T-junctions,
which increases the flexibility of the geometric model significantly.
In this IETI-DP variant, the construction of the coarse space is based on the idea of ``fat vertices'': We consider every basis function that is supported on a vertex or T-junction as primal degree of freedom. The numerical experiments indicate that a similar condition number bound to~\eqref{estimate} might hold.  

The remainder of this paper is organized as follows. In Secion~\ref{sec:problem setting} we describe the model problem. In Section~\ref{sec:solver} we introduce the IETI-DP solver and we end this paper with numerical experiments in Section~\ref{sec:numerics}.

\section{The problem setting}
\label{sec:problem setting}
Let $\Omega \subset \mathbb{R}^2$ be open, simply connected and bounded with Lipschitz boundary $\partial \Omega$. $L_2(\Omega)$ and $H^1(\Omega)$
are the common Lebesgue and Sobolev spaces. As usual, $H^1_0(\Omega) \subset H^1(\Omega)$ denotes the subspace of functions that vanish on $\partial \Omega$.

We consider the following model problem:
Find $u \in H^1_0(\Omega)$ such that 
\begin{align}
	\label{continousProb}
	\int_{\Omega}^{} \nabla u \cdot \nabla v \; \mathrm{d}x = \int_{\Omega}^{} fv \; \mathrm{d}x \qquad\text{for all}\qquad v \in H^1_0(\Omega)
\end{align}
with a given source function $f \in L_2(\Omega)$.  
We assume that $\Omega$ is a composition of $K$ non-overlapping patches $\Omega^{(k)}$, where 
every patch $\Omega^{(k)}$ is parametrized
by a 
geometry function
\begin{align}
	G_k:\widehat{\Omega}:=(0,1)^2 \rightarrow \Omega^{(k)}:=G_k(\widehat{\Omega}) \subset \mathbb{R}^2, 
\end{align}
that has a continuous extension to the closure of $\widehat\Omega$ and
such that
$\nabla G_k \in L_\infty(\widehat{\Omega})$
and
$(\nabla G_k)^{-1} \in L_\infty(\widehat{\Omega})$.

We consider the case where the pre-images of the (Dirichlet) boundary consist of whole edges.
The indices of neighboring patches $\Omega^{(\ell)}$ of $\Omega^{(k)}$, that share at least a part of their boundaries, is collected in the set
\[
\mathcal N_\Gamma(k):=
\{ \ell \neq k \;:\; \mbox{meas }({ \partial \Omega^{(k)} } \cap { \partial \Omega^{(\ell)} } ) > 0 \},
\]
where $\mbox{meas }T$ is the measure of $T$.
For any $\ell \in \mathcal{N}_\Gamma(k)$, we
write  $\Gamma^{(k,\ell)}={\partial \Omega^{(k)}} \cap {\partial \Omega^{(\ell)}}$.
The endpoints of ${\partial \Omega^{(k)}} \cap {\partial \Omega^{(\ell)}}$
that are not located on the (Dirichlet) boundary of $\Omega$
are referred to as junctions. A junction could be a common vertex or
a T-junction.

For the IgA discretization spaces, we first construct a B-spline space $\widehat{V}^{(k)}$ on the parameter domain $\widehat{\Omega}$ by tensorization of two univariate B-spline spaces. 
The function spaces on the physical domain are then defined by the
pull-back principle: $V^{(k)} := \widehat{V}^{(k)} \circ G_k^{-1}$.

The product of the local spaces gives the global approximation space 
$
	V := V^{(1)} \times \dots \times V^{(K)}.
$
On this discretization space, we could introduce the SIPG formulation, cf.~\cite{Arnold:1982, SchneckenleitnerTakacs:2020}. 
Since we are interested in a domain decomposition approach, we need
patch-local formulations of SIPG.
\section{The dG IETI-DP solver}
\label{sec:solver}
For those patch-local formulations, we
adapt the ideas of~\cite{DryjaGalvis:2013,HoferLanger:2016a,Hofer:2016a} and others.
We choose local function spaces $V_{e}^{(k)}$ to be the product space of $V^{(k)}$ and the neighboring trace spaces~$V^{(k,\ell)}$, which are
the restrictions of $V^{(\ell)}$ to $\Gamma^{(k,\ell)}$.
A function $v_e^{(k)} \in V_e^{(k)}$ is represented as a tuple
$
	\label{def:representation}
	v_e^{(k)} = 
	\left(
	v^{(k)}, (v^{(k,\ell)})_{\ell\in \mathcal{N}_\Gamma(k)}
	\right),
$
where $v^{(k)} \in V^{(k)}$ and $v^{(k,\ell)} \in V^{(k,\ell)}$. Note that the traces of the basis functions for $V^{(\ell)}$ restricted to $\Gamma^{(k,\ell)}$ form a basis of $V^{(k,\ell)}$.
The basis for $V_e^{(k)}$ consists of the basis functions of $V^{(k)}$ and the basis functions for $V^{(k,\ell)}$. The basis functions on $V^{(k,\ell)}$ are usually visualized as living on artificial interfaces.

\newpage

On each patch, we consider the local problem:
Find $u_e^{(k)} \in V_e^{(k)}$ such that 
\begin{align*}
	a_e^{(k)}(u_e^{(k)},v_e^{(k)})=\langle f_e^{(k)},v_e^{(k)}\rangle
	\quad \mbox{for all}\quad v_e^{(k)} \in V_e^{(k)},
\end{align*}
where
\[
\begin{aligned}
	a_e^{(k)}(u_e^{(k)},v_e^{(k)}) &:= a^{(k)}(u_e^{(k)},v_e^{(k)}) + m^{(k)}(u_e^{(k)},v_e^{(k)}) + r^{(k)}(u_e^{(k)},v_e^{(k)}), \\
	\langle f_e^{(k)},v_e^{(k)}\rangle & := \int_{\Omega^{(k)}} f v^{(k)} \mathrm dx,\\
	a^{(k)}(u_e^{(k)},v_e^{(k)}) &:= \int_{\Omega^{(k)}} \nabla u^{(k)} \cdot \nabla v^{(k)} \; \textrm{d}x, \\
	m^{(k)}(u_e^{(k)},v_e^{(k)}) &:= \sum_{\ell \in \mathcal{N}_\Gamma(k)} \int_{\Gamma^{(k,\ell)}}  
	\frac{\partial u^{(k)}}{\partial n_k}(v^{(k,\ell)} - v^{(k)}) 
	\; \textrm{d}s, \\
	&\qquad + \sum_{\ell \in \mathcal{N}_\Gamma(k)} \int_{\Gamma^{(k,\ell)}}  
	\frac{\partial v^{(k)}}{\partial n_k}(u^{(k,\ell)} - u^{(k)})
	 \; \textrm{d}s, \\
	r^{(k)}(u_e^{(k)},v_e^{(k)}) &:= \sum_{\ell \in \mathcal{N}_\Gamma(k)} \int_{\Gamma^{(k,\ell)}} \frac{\delta p^2}{h_{k\ell}}   
	(u^{(k,\ell)} - u^{(k)})(v^{(k,\ell)} - v^{(k)}) \; \textrm{d}s
\end{aligned}
\]
and $n_k$ denotes the outward unit normal vector and $\delta$ is
the dG penalty parameter, which has to be chosen large enough in
order to guarantee that the bilinear form
$a_e^{(k)}(\cdot,\cdot)$ is coercive. In~\cite{Takacs:2019b}, it
was shown that $\delta$ can be chosen independently of $p$.

The discretization of $a_e^{(k)}(\cdot,\cdot)$ and $\langle f_e^{(k)},\cdot \rangle$
gives a local system, which we write as 
\begin{equation}\label{linsys:local}
\left(
\begin{array}{cccc}
	A_{\mathrm I \mathrm I}^{(k)} & A_{\mathrm I \mathrm \Gamma}^{(k)}\\
	A_{\mathrm \Gamma \mathrm I}^{(k)} & A_{\mathrm \Gamma \mathrm \Gamma}^{(k)}\\
\end{array}
\right)
\left(
\begin{array}{c}
	\underline{u}_\mathrm{I}^{(k)} \\ \underline{u}_{\mathrm \Gamma}^{(k)}
\end{array}
\right)
=
\left(
\begin{array}{c}
	\underline{f}_{\mathrm I}^{(k)}\\\underline{f}_{\mathrm \Gamma}^{(k)}
\end{array}
\right),
\end{equation}
where the index $\mathrm{I}$ refers to the basis functions that are
only supported in the interior of $\Omega^{(k)}$ and the index
$\mathrm \Gamma$ refers to the remaining basis functions, i.e., those living
on the patch boundary and on the artificial interfaces.
We eliminate the interior degrees of freedom in \eqref{linsys:local} for every $k = 1,\dots, K$ to get the block diagonal Schur complement system
\begin{equation}\label{eq:schursys}
	S \underline w = \underline g,
\end{equation}
where the individual blocks of $S$ are given by
$S^{(k)} = A_{\mathrm \Gamma \mathrm \Gamma}^{(k)} -
A_{\mathrm \Gamma \mathrm I}^{(k)}
\big(A_{\mathrm I \mathrm I}^{(k)}\big)^{-1}
A_{\mathrm I \mathrm \Gamma}^{(k)}$.

The IETI-DP method requires carefully selected primal degrees of freedom to be solvable. We choose the degrees of freedom associated to the basis functions which are non-zero on a junction to be primal.
For every standard corner, we only have one primal degree of freedom
per patch, as in~\cite{SchneckenleitnerTakacs:2020}. On a T-junction however, the
number of non-zero basis functions grows linearly with $p$. Since
we take all of them, we refer to ``fat vertices'' in this context.

$C = \mbox{diag }(C^{(1)}, \dots, C^{(K)})$ is the constraint matrix,
i.e., it is defined such that $C \underline w = 0$ if and only if
the associated function $w$ vanishes at the primal degrees of freedom.
The matrix $\Psi$ represents the energy minimizing basis functions
for the space of primal degrees of freedom.

\begin{SCfigure}[1][htb]
	\scalebox{0.75}{
	\begin{tikzpicture}
		\fill[gray!20] (-0.2,0) -- (4.7,0) -- (4.7,1.2) -- (-0.2,1.2);
		\fill[gray!20] (-0.2,-1.5) -- (1.5,-1.5) -- (1.5,-3.2) -- (-0.2,-3.2);
		\fill[gray!20] (4.7,-1.5) -- (3.0,-1.5) -- (3.0,-3.2) -- (4.7,-3.2);
		
		\draw (-0.2,0) -- (4.7,0) node at (2.5,0.55) {$\Omega^{(1)}$};
		\draw (-0.2,-1.5) -- (1.5, -1.5) -- (1.5,-3.2) node at (0.7,-2.35) {$\Omega^{(2)}$};
		\draw (4.7,-1.5) -- (3.0,-1.5) -- (3.0,-3.2) node at (3.9,-2.35) {$\Omega^{(3)}$};
		
		\draw (-0.2,-0.4) -- (1.5,-0.4) {};			
		\draw (3.0,-0.4) -- (4.7,-0.4) {};
		\draw (-0.2,-1.1) -- (2.0,-1.1) {};
		\draw (1.9,-1.5) -- (1.9,-3.2) {};
		\draw (2.5,-1.1) -- (4.7,-1.1) {};
		\draw (2.6,-1.5) -- (2.6,-3.2) {};
		
		\draw (0.0,0.0) node[circle, fill, inner sep = 2.5pt] (A1) {};
		\draw (1.5,0) node[circle, fill, inner sep = 2.5pt] (A2) {};
		\draw (3.0,0.0) node[circle, fill, inner sep = 2.5pt] (A4) {};
		\draw (4.5,0.0) node[circle, fill, inner sep = 2.5pt] (A5) {};
		\draw (3.0,-0.4) node[star, fill, inner sep = 2pt] (A6) {};
		\draw (3.75,-0.4) node[star, fill, inner sep = 2pt] (A7) {};
		\draw (4.5,-0.4) node[star, fill, inner sep = 2pt] (A8) {};
		\draw (0.75,-0.4) node[diamond, fill, inner sep = 2pt] (A9) {};
		\draw (1.5,-0.4) node[diamond, fill, inner sep = 2pt] (A10) {};

		\draw (1.5,-1.5) node[diamond, fill, inner sep = 2pt] (B1) {};
		\draw (1.5,-2.2) node[diamond, fill, inner sep = 2pt] (B2) {};
		\draw (1.5,-3.0) node[diamond, fill, inner sep = 2pt] (B3) {};
		\draw (0.75,-1.5) node[diamond, fill, inner sep = 2pt] (B4) {};
		\draw (1.9,-1.5) node[star, fill, inner sep = 2pt] (B5) {};
		\draw (1.9,-2.0) node[star, fill, inner sep = 2pt] (B6) {};
		\draw (1.9,-2.75) node[star, fill, inner sep = 2pt] (B7) {};
		\draw (0.0,-1.1) node[circle, fill, inner sep = 2.5pt] (B8) {};
		\draw (1.0,-1.1) node[circle, fill, inner sep = 2.5pt] (B9) {};
		\draw (2.,-1.1) node[circle, fill, inner sep = 2.5pt] (B10) {};

		\draw (3.0,-1.5) node[star, fill, inner sep = 2pt] (C1) {};
		\draw (3.0,-2.0) node[star, fill, inner sep = 2pt] (C2) {};
		\draw (3.0,-2.75) node[star, fill, inner sep = 2pt] (C3) {};
		\draw (3.75,-1.5) node[star, fill, inner sep = 2pt] (C4) {};
		\draw (4.5,-1.5) node[star, fill, inner sep = 2pt] (C5) {};
		\draw (2.6,-1.5) node[diamond, fill, inner sep = 2pt] (C6) {};
		\draw (2.6,-2.2) node[diamond, fill, inner sep = 2pt] (C7) {};
		\draw (2.6,-3.0) node[diamond, fill, inner sep = 2pt] (C8) {};
		\draw (2.5,-1.1) node[circle, fill, inner sep = 2.5pt] (C9) {};
		\draw (3.5,-1.1) node[circle, fill, inner sep = 2.5pt] (C10) {};
		\draw (4.5,-1.1) node[circle, fill, inner sep = 2.5pt] (C11) {};
		
		\draw[dotted, <->, line width = 1pt, latex-latex, bend left]
		(A1) edge (B8) 
		(A7) edge (C4) (A8) edge (C5)
		(B6) edge (C2) (B7) edge (C3); 
		
		\draw[dotted, <->, line width = 1pt, latex-latex, bend right]
		(A9) edge (B4) 
		(A5) edge (C11)
		(B2) edge (C7) (B3) edge (C8);
		
		\draw[<->, line width = 1pt, latex-latex, bend right]
		(A2) edge (B9) (A4) edge (B10) (A10) edge (B1) (B1) edge (C6) (C1) edge (B5);
		
		\draw[<->, line width = 1pt, latex-latex, bend left]
		(A2) edge (C9) (A4) edge (C10) (A6) edge (C1);
	\end{tikzpicture}
}
	\caption{Action of matrix $B$ (dotted lines) \\ 
	and primal degrees of freedom (solid lines)} \label{fig:ommiting}
\end{SCfigure}
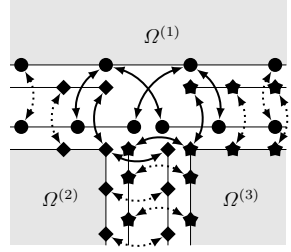

Furthermore, we introduce the jump matrix $B$, which models the jumps
of the functions between the patch boundaries and the the associated
artificial interfaces. Each row corresponds to one degree of freedom
(coefficient for a basis function) on the the patch boundary and one
artificial interface; as usual, each row has only two non-zero
coefficients that are $-1$ and $1$.
Primal degrees of freedom are excluded. For
a visualization, see Fig.~\ref{fig:ommiting}, where the primal degrees
of freedom are marked with solid lines
and the dotted arrows show the action of the jump matrix $B$.
The basis functions on the artificial interfaces are labeled with the
same symbols from the original spaces.

The following problem is equivalent to the SIPG discretization of~\eqref{continousProb}, cf. \cite{MandelDohrmannTezaur:2005a}: Find
$(\underline w_\Delta, \underline \mu, \underline w_\Pi, \underline \lambda)$ such that
\[
\begin{pmatrix}
	S & C^\top &                  & B^\top           \\
	C &        &                  &                  \\
	&          & \Psi^\top S \Psi\;\;\; & (B\Psi)^\top \\
	B &        & B \Psi           &                  \\
\end{pmatrix}
\begin{pmatrix}
	\underline{w}_\Delta  \\
	\underline{\mu}  \\
	\underline{w}_\Pi  \\
	\underline{\lambda}  \\
\end{pmatrix}
=
\begin{pmatrix}
	\underline{g}  \\
	0  \\
	\Psi^\top \underline{g}  \\
	0  \\
\end{pmatrix}.
\]
We obtain the solution of the original problem by
$\underline{w} = \underline{w}_\Delta + \Psi \underline{w}_\Pi$.
We build a Schur complement of this system to get the linear problem   
\begin{equation}
	\label{IETIProblem}
	F \; \underline{\lambda} = \underline{d}.
\end{equation}
We solve~\eqref{IETIProblem} with a preconditioned conjugate gradient (PCG) solver with the scaled Dirichlet preconditioner 
\[
M_{\mathrm{sD}} := B D^{-1} S  D^{-1} B^\top,
\]
where $D$ is a diagonal matrix
defined based on the principle of multiplicity scaling, cf.~\cite{SchneckenleitnerTakacs:2019, Pechstein:2013a}.

\section{Numerical results}
\label{sec:numerics}
\begin{figure}[h]
	\centering
	\begin{subfigure}{.27\textwidth}
		\centering
		\includegraphics[height=3cm]{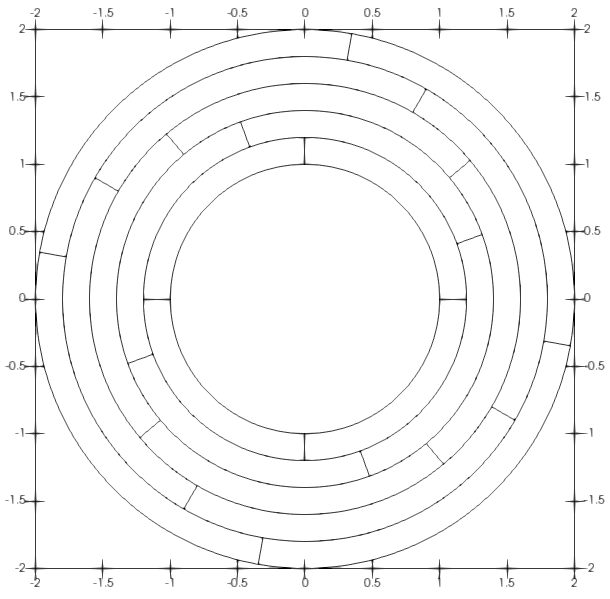}
		{\tiny \caption{Ring} \label{subfig:Ring}}
	\end{subfigure} 
	\hspace{1.5cm}
	\begin{subfigure}{.27\textwidth}
		\centering
		\includegraphics[height=3cm]{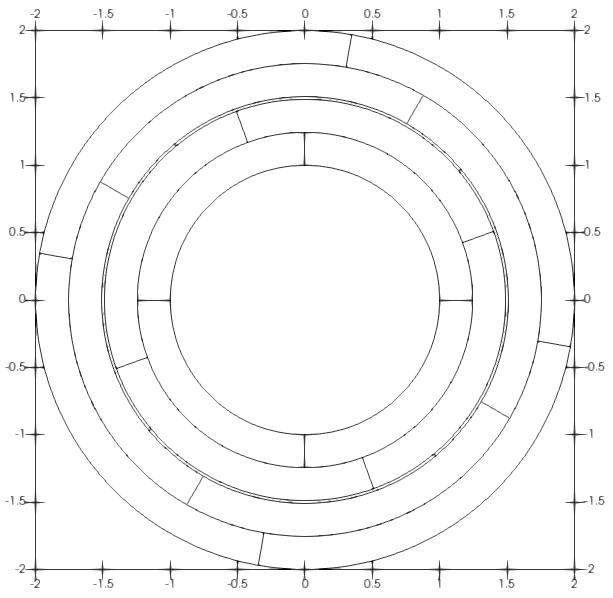}
		{\tiny \caption{Ring with thin gap}
			\label{subfig:Ring thin}}
	\end{subfigure} 
	\caption{Computational domains and the decomposition into patches}
	\label{fig:computational domains}
\end{figure}

We consider the model problem 
\[
\begin{aligned}
	- \Delta u(x,y) & = 2\pi^2 \sin(\pi x)\sin(\pi y) &&\qquad \mbox{for}\quad (x,y)\in\Omega \\
	u & = 0 &&\qquad \mbox{on}\quad \partial\Omega,
\end{aligned}
\]
on the geometries depicted in Fig.~\ref{fig:computational domains}.
Both represent the same computational domain with an inner radius of $1$ and an outer radius of $2$. The ring in Fig.~\ref{subfig:Ring} consists of $20$ patches each of which has a width of $0.2$. For the ring in Fig.~\ref{subfig:Ring thin}, the thin layer of patches has a width of $0.02$, the other
layers have a correspondingly larger width. 
We use NURBS of degree $2$ to parametrize all patches. In the coarsest setting, i.e., $r = 0$, the discretization spaces on all patches consist of global polynomials only. The discretization spaces for $r = 1,2,3,\dots$ are obtained by uniform refinement steps. We use a PCG solver to solve system~\eqref{IETIProblem} with the preconditioner $M_{\mathrm{sD}}$ and to estimate the condition number $\kappa(M_\mathrm{sD} F)$, where we use the zero vector as initial guess.
All experiments are carried out in the C++ library G+Smo, cf.~\cite{gismoweb} and are executed on the Radon1\footnote{https://www.ricam.oeaw.ac.at/hpc/} cluster in Linz. 

In the Table~\ref{tab:Ring},
we report on the iteration counts (it) and the condition numbers $(\kappa)$
for various refinement levels $r$ and various spline degrees $p$,
where we chose $C^s$-smoothness with $s=p-1$ within the patches. The tables show the expected behavior with respect to $h$. The condition number decreases when we increase the spline degree $p$, which is better than one would expect from the theory in~\cite{SchneckenleitnerTakacs:2020}.
Although the width of the thin patches in Fig.~\ref{subfig:Ring thin} is one tenth of the width of the patches in Fig.~\ref{subfig:Ring}, the condition number grows only by a factor between $5$ and $6$. 
Also the iteration counts grow only mildly.
%

\begin{table}[ht]
	\newcolumntype{L}[1]{>{\raggedleft\arraybackslash\hspace{-1em}}m{#1}}
	\centering
	\renewcommand{\arraystretch}{1.25}
		\begin{tabular}{p{.7em}||L{1.3em}L{2.3em}|L{1.3em}L{2.3em}|L{1.3em}L{2.3em}|L{1.3em}L{2.3em}||
		                  L{1.3em}L{2.3em}|L{1.3em}L{2.3em}|L{1.3em}L{2.3em}|L{1.3em}L{2.3em}}
			\toprule
			& \multicolumn{8}{c||}{Fig.~\ref{subfig:Ring}}
			& \multicolumn{8}{c}{Fig.~\ref{subfig:Ring thin}}
			\\
			& \multicolumn{2}{c|}{$p=2$}
			& \multicolumn{2}{c|}{$p=3$}
			& \multicolumn{2}{c|}{$p=6$}
			& \multicolumn{2}{c||}{$p=7$}
			& \multicolumn{2}{c|}{$p=2$}
			& \multicolumn{2}{c|}{$p=3$}
			& \multicolumn{2}{c|}{$p=6$}
			& \multicolumn{2}{c}{$p=7$} \\
			$r$
			& it & $\kappa$
			& it & $\kappa$
			& it & $\kappa$
			& it & $\kappa$
			& it & $\kappa$
			& it & $\kappa$
			& it & $\kappa$
			& it & $\kappa$  \\
			\midrule
			$4$  & $9$ & $3.7$ & $9$ & $3.5$ 
			& $8$ & $2.4$ & $8$ & $2.1$ 
			& $12$ & $18.0$ & $13$ & $18.0$ 
		  & $13$ & $12.9$ & $11$ & $11.6$
			\\
			$5$  & $10$ & $4.6$ & $10$ & $4.5$ 
			& $9$ & $3.8$ & $9$ & $3.5$ 
			& $20$ & $24.1$ & $19$ & $23.4$ 
			& $19$ & $19.3$ & $18$ & $18.0$
			\\
			$6$  & $10$ & $5.8$ & $10$ & $5.5$ 
			& $10$ & $4.9$ & $10$ & $4.8$ 
			& $22$ & $31.7$ & $22$ & $29.9$ 
			& $21$ & $25.9$ & $20$ & $24.8$
			\\
			$7$  & $11$ & $6.3$ & $11$ & $6.2$ 
			& $10$ & $5.6$ & $10$ & $5.5$ 
			& $24$ & $37.2$ & $24$ & $36.3$ 
			& $22$ & $31.3$ & $22$ & $30.1$
			\\
			$8$  & $11$ & $6.7$ & $11$ & $6.7$ 
			& $11$ & $6.3$ & $10$ & $5.6$ 
			& $24$ & $43.2$ & $24$ & $42.3$ 
			& $24$ & $36.5$ & $24$ & $31.6$
			\\
			\bottomrule
	\end{tabular}
	\caption{Iterations (it) and condition numbers ($\kappa$)
		\label{tab:Ring}}
\end{table}

The Table~\ref{tab:Ring time} presents the parallel solving times for $n$ processors. We only consider
the domain in Fig.~\ref{subfig:Ring} again with $s=p-1$. We see that the speedup rate with respect to $n$ is a bit smaller than the expected rate of $2$. This is probably caused by the rather small number of patches in the computational domain. 

	\begin{table}[ht]
	\newcolumntype{L}[1]{>{\raggedleft\arraybackslash\hspace{-1em}}m{#1}}
	\centering
	\renewcommand{\arraystretch}{1.25}
		\begin{tabular}{l||L{2.5em}|L{2.5em}|L{2.5em}|L{2.5em}|L{2.5em}||
		                  L{2.5em}|L{2.5em}|L{2.5em}|L{2.5em}|L{2.5em}}
			\toprule
			& \multicolumn{5}{c||}{$p=3$}
			& \multicolumn{5}{c}{$p=7$} \\
			$r$
			& \multicolumn{1}{c|}{$n=1$}
			& \multicolumn{1}{c|}{$n=2$}
			& \multicolumn{1}{c|}{$n=4$}
			& \multicolumn{1}{c|}{$n=8$}
			& \multicolumn{1}{c||}{$n=16$}
			& \multicolumn{1}{c|}{$n=1$}
			& \multicolumn{1}{c|}{$n=2$}
			& \multicolumn{1}{c|}{$n=4$}
			& \multicolumn{1}{c|}{$n=8$}
			& \multicolumn{1}{c}{$n=16$}  \\
			\midrule
			$6$  & $3.8$ & $2.8$ & $2.4$ & $ 1.2$ & $0.8$
			     & $10.0$ & $6.5$ & $5.0$ & $2.5$ & $1.75$ \\
			$7$  & $24.0$ & $16.1$ & $13.6$ & $6.4$ & $4.1$ 
			     & $47.0$ & $31.7$ & $26.1$ & $12.2$ & $9.3$\\
			$8$  & $107.0$ & $81.4$ & $66.8$ & $29.5$ & $19.5$ 
			     & $220.0$ & $158.7$ & $129.4$ & $56.7$ & $45.4$\\
			\bottomrule
	\end{tabular}
	\caption{Solving times (sec.); Fig.~\ref{subfig:Ring}
		\label{tab:Ring time}}
	\end{table}

In Table~\ref{tab:Ring cont} we report on the iteration counts and the condition numbers for the decomposition in Fig.~\ref{subfig:Ring} when we change the smoothness $s$ of the B-splines within the patches. The numbers in the table show the behavior for $r=5$. We see that for a fixed smoothness $s$ the condition number grows slightly with respect to the spline degree $p$. For a fixed degree $p$, we observe a decline in the condition number when we increase the smoothness~$s$. 

\begin{table}[ht]
	\newcolumntype{L}[1]{>{\raggedleft\arraybackslash\hspace{-1em}}m{#1}}
	\centering
	\renewcommand{\arraystretch}{1.25}
		\begin{tabular}{l||L{1.8em}L{2.em}|L{1.8em}L{2.em}|L{1.8em}L{2.em}|L{1.8em}L{2.em}|L{1.8em}L{2em}|L{1.8em}L{2em}}
			\toprule
			& \multicolumn{2}{c|}{$p=2$}
			& \multicolumn{2}{c|}{$p=3$}
			& \multicolumn{2}{c|}{$p=4$}
			& \multicolumn{2}{c|}{$p=5$}
			& \multicolumn{2}{c|}{$p=6$}
			& \multicolumn{2}{c}{$p=7$} \\
			$s$ 
			& it & $\kappa$
			& it & $\kappa$
			& it & $\kappa$
			& it & $\kappa$
			& it & $\kappa$
			& it & $\kappa$ \\
			\midrule
			$0$  & $10$ & $5.0$ & $10$ & $5.3$ & $10$ & $5.4$ & $10$ & $5.5$ & $10$ & $5.6$ & $10$ & $5.6$ \\
			$1$  & $10$ & $4.6$ & $10$ & $5.2$ & $10$ & $5.3$ & $10$ & $5.4$ & $10$ & $5.5$ & $10$ & $5.5$ \\
			$2$  & $ $ & $ $ & $10$ & $4.5$ & $10$ & $5.0$ & $10$ & $5.3$ & $10$ & $5.4$ & $10$ & $5.5$ \\
			$3$  & $ $ & $ $ & $ $ & $ $ & $9$ & $4.2$ & $10$ & $4.9$ & $10$ & $5.1$ & $10$ & $5.3$ \\
			$4$  & $ $ & $ $ & $ $ & $ $ & $ $ & $ $ & $9$ & $4.0$ & $10$ & $4.7$ & $10$ & $5.0$ \\
			$5$  & $ $ & $ $ & $ $ & $ $ & $ $ & $ $ & $ $ & $ $ & $9$ & $3.8$ & $9$ & $4.5$ \\
			$6$  & $ $ & $ $ & $ $ & $ $ & $ $ & $ $ & $ $ & $ $ & $ $ & $ $ & $9$ & $3.5$ \\
			\bottomrule
	\end{tabular}
	\caption{Iterations (it) and condition number ($\kappa$); $r = 5$; Fig.~\ref{subfig:Ring}
		\label{tab:Ring cont}}
\end{table}


\begin{acknowledgement}
	The first author was supported by the Austrian Science Fund (FWF): S117 and W1214-04. Also, the second author has received support from the Austrian Science	Fund (FWF): P31048.
\end{acknowledgement}

%
%
\bibliography{literature}

\end{document}